# A Performance-Based Framework for Bridge Preservation Based on Damage-Integrated System-Level Behavior

TRB Paper # 14-2715

Date Submitted: November 12, 2013

| | |
|---:|---:|
| Abstract: | 244 |
| Main Text: | 5534 |
| Tables | 250 |
| Figures: | 2000 |
| Total Word Count: | 8028 |

**Authors:**

**Amir Gheitasi**, Graduate Research Assistant, University of Virginia
Department of Civil and Environmental Engineering
351 McCormick Road, Charlottesville, VA, 22904-4742, USA
Phone: (906) 370-4557
Fax: (434) 982-2951
Email: agheitasi@virginia.edu (Corresponding Author)

**Devin K. Harris**, Assistant Professor, University of Virginia
Department of Civil and Environmental Engineering
351 McCormick Road, Charlottesville, VA, 22904-4742, USA
Phone: (434) 924-6373
Fax: (434) 982-2951
Email: dharris@virginia.edu

**ABSTRACT**

The safety and condition of transportation infrastructure has been at the forefront of national debates in recent times due to catastrophic bridge failures, but the issue has been a longstanding challenge for transportation agencies for many years as resources continue to diminish. The performance of this infrastructure has a direct influence on the lives of most of citizens in developed regions by providing a critical lifeline between communities and the transportation of goods and services, and as a critical component of the transportation network, bridges have received a lot of attention regarding condition assessment and maintenance practices. Despite successful implementation of advanced evaluation techniques, what is still lacking is a fundamental understanding of the system behavior in the presence of deteriorating conditions that can be used for preservation decision-making.

   This paper aims to present a performance-based framework that can be used to characterize the behavior of in-service bridge superstructures. In order to measure the bridge system performance with deteriorating conditions, system-level behavior of a representative composite steel girder bridge, degraded with three common damage scenarios was investigated in this study. Results obtained from validated numerical analysis demonstrated significant impact of integrated damage mechanisms on the ultimate capacity, redundancy and system ductility of the simulated bridge superstructure. It is expected that the proposed framework for evaluating system behavior will provide a first step for establishing a critical linkage between design, maintenance, and rehabilitation of highway bridges, which are uncoupled in current infrastructure decision-making processes.

**STATE OF PRACTICE**
An efficient and well maintained transportation infrastructure is an essential component to the economic health of the United States not only by providing a corridor for the transportation of goods and people (1, 2), but also serves a coast to coast and border to border passageway for the nation's military (3). Historically the highway system, which includes both roads and bridges, has flown underneath the radar of the public opinion when compared to other engineering marvels, but is easily recognized as a critical component in the everyday lives of most people. Tragic events such as the I-35W bridge collapse in August 2007 (4) and the more recent I-5 collapse in may 2013 (5) have brought the challenges associated with an aging infrastructure to the forefront of public scrutiny. The cause of this scrutiny can be attributed to the fact that these failures often result in loss of human life and significant economic hardship to the surrounding communities (6).

The occurrence of bridge failures is somewhat rare, sometimes even preventable with proper maintenance and repair, especially when coupled with reliable mechanisms for damage diagnosis. These uncommon failures can be attributed to unforeseen events such as impact, fires or flooding and are typically not attributed exclusively to deterioration (7, 8). Nevertheless, it is the condition states associated with deterioration that represents the greatest challenges for transportation agencies across the country. With a national inventory of more than 600,000 bridges (9), 25 percent of which are classified as structurally deficient (11%) or functionally obsolete (14%), strategies and resources for maintenance are a growing challenge for federal, state and local governments, especially considering that many bridges are reaching or exceeding their design service lives of 50 years. While it is not feasible to immediately repair all of the deficient bridges, this deteriorating condition does underscore the importance of quality inspection and performance assessment mechanisms to prioritize these repair efforts.

**COMMON DETERIORATIONS IN HIGHWAY BRIDGES**
Within a given genre of the bridge superstructures, there are certain similarities in materials, geometry, and configuration, which make them serve common functionalities. Under the premise of a rational structural design, the service lives of these bridges are governed by the operating environment, load effects and history, and maintenance and preservation practices; with really only the last factor being under the owner's influence. Successful maintenance practices require knowledge of condition state of the bridge system and also an understanding of the impacts of all possible damage and deterioration mechanisms on the overall system performance and serviceability. To some extent, this is accomplished through traditional inspection and load rating practices employed by most transportation agencies, but current inspection methods have limitations in their ability to detect damage while load rating methods follow the element-level analysis rationale used in design which does not truly account for the complex system-level behavior.

The type of damage and deterioration conditions associated with each bridge system varies from one structural type to another (10). In composite steel girder bridges for instance, much of the degradation often manifests in the steel girders as corrosion, section loss, and fatigue cracking. Moreover, damage mechanisms such as corroded reinforcing steel, subsurface delaminations, spalled sections, and scaling are the most well-known sources of degradation associated with the concrete decks. Other mechanisms such as frozen bearings, over-height vehicular impacts, and bridge settlement, which address the overall system behavior, are also common in this type of bridge superstructure systems.

**PROBLEM TO ADDRESS**

Over the past few years, a number of research efforts have concentrated on developing an applicable mechanism to qualitatively and quantitatively characterize the main sources of defects associated with the transportation infrastructures, including bridges. The concept of structural health monitoring is one of these developed methodologies that have come to the forefront of the research community with a primary focus of accurately monitoring in-situ behavior to assess in-service performance, detecting damage, and determining condition of a structure. Recent advances in non-destructive evaluation (NDE) (11, 12) along with novel technologies such as fiber optic sensors, distributed wireless sensors and networks (13-15), have furthered the science of assessment, allowing for more accurate quantification of visible deterioration mechanisms and improved confidence in locating internal deterioration mechanisms.

With the sheer volume of in-service bridges in the United States, the main challenges with the current practices of condition assessment are the potential for large amounts of collected data and the manpower and skills required to interpret and manage them. Moreover, the basic question that still needs to be answered is how this collected data can be used to correlate the impact of existing damage and deteriorating conditions on the performance of highway bridges. In fact, what transportation decision makers are lacking is a fundamental understanding of system-level behavior characteristics of the deteriorated bridges during their varying life stages, and the correlation of this behavior to safety, remaining life, and preservation strategies that are essential to defining a bridge's life-cycle.

The focus of this paper is to provide a mechanism for understanding the complex system-level interaction inherent to the bridges and characterizing the influence of damage and deteriorating conditions on these interactions as well as overall system behavior and performance of in-service bridge superstructures. This study has been limited to composite steel girder bridges, as they represent one of the most common structure types in service (9), but the proposed approach is generic, allowing for extrapolation across other bridge types in the inventory having a wide variety of structure types, component materials, existing conditions, and operational environments. This study proposes a first step in a philosophy change aimed at establishing damage-integrated baseline system performance measures, from which preservation decisions can be made for in-service bridges. In current preservation practices, there exists a lack of baseline behavior characterization for in-service bridges beyond the component-based behavior assumed in design, especially when the damage is present.

**INVESTIGATION APPROACH**

With the goal of establishing a performance-based framework to evaluate the behavior of in-service bridge superstructures, the investigation approach of this study has been categorized into four different phases to accomplish this goal. Figure 1 illustrates the schematic of the proposed generic framework that can be applied to any type of bridge superstructure. With the limited data available on the ultimate capacity testing of in-service bridges with damage conditions and the lack of understanding regarding the system-level behavior of bridge superstructures, this framework aims to describe a rationale for linking these components together to yield a rational system-level model with integrated damage mechanisms.

As illustrated, the first phase focuses on the development of undamaged element-level numerical models representing the main structural components of the desired bridge system. Extensive experimental data in the literature are available to validate the accuracy and consistency of the developed models. This phase aims to characterize the behavior and failure modes of the intact bridge sub-components, with consideration of both material and geometric non-linearities.

Following validation of the intact models within the element domain, phase two focuses on the development of system-level numerical models representing ideal intact bridge superstructures. These models can be loaded to their ultimate capacity to define their full non-linear system behavior. Using the non-linear system behavior highlights the critical behavioral stages inherent to a particular bridge system and also allows for correlation with the expected element-level response. Ideally complete experimental data of newly constructed bridges tested to failure is needed to validate, update, and/or refine the established system-level models; however, the primary challenge associated with model validation of bridges is the lack of complete datasets, especially those for full non-linear behavior characterization (16-18). The limited pool of data exists in the literature can be used in this study to complete the skeleton of the proposed framework, which could also be supplemented with focused laboratory or field investigations to fill the void of available experimental data. The outcome of this phase will be a mechanism to describe the behavior of bridge systems based on their inherent level of system redundancy that are expected to be unique to each bridge superstructure type and design characteristics. This component will be essential to the latter stages of the investigation as it provides a baseline for as-designed behavior from which actual in-service behavior can be referenced to define the influence of damage and deterioration (phase IV).

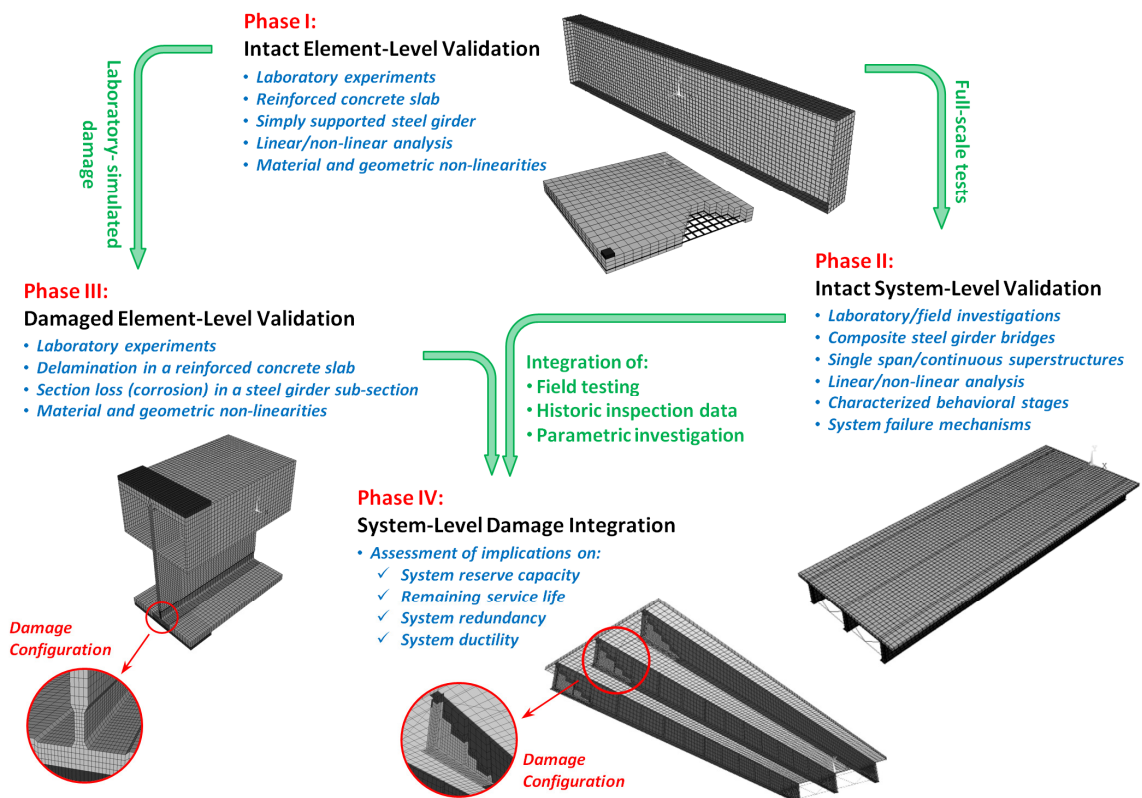

**FIGURE 1 Schematic representation of the proposed framework.**

Following the first two phases, the purpose of the third phase of the proposed framework is to establish an effective constitutive behavior of the individual bridge components that are affected by common damage and deteriorating mechanisms. Despite the variety of sources that may cause each damage scenario, from a mechanics perspective it is their influence on the structural behavior that is of primary concern. Critical to this phase is the strategy to leverage modeling techniques appropriate for each damage mechanism considered. Validation of the damage-integrated elements through available experimental data is essential to the last phase of this investigation. Upon completion of the damaged element-level validation, the damage modeling strategies will be integrated into the system-level models to investigate their influence on system behavior (Phase IV). Parametric investigations can then be performed to quantify the influence of damage mechanisms on the system-level behavioral stages defined in Phase II and establish a measure of remaining life and susceptibility to failure.

In order to demonstrate the applicability of the proposed framework, the associated methodology has been applied in this study to the representative category of composite steel girder bridges. Numerical models of the intact and damaged structural systems within both element-level and system-level domains have been developed based on the existing literature. The commercial finite element (FE) computer package, ANSYS (19), was used to create and refine each of these numerical models. The accuracy and validity of the FE simulation and analysis were also investigated through a comparison of the numerical results to the corresponding experimental data. The results presented in this paper are a representative subset of the analyses performed within the proposed framework.

For other types of bridge superstructures, the proposed approach can be updated for specific structural components and system configurations to study the effect of deteriorating mechanisms on the performance of the bridge system. These damage-integrated system models represent a critical component to a performance-based preservation framework, as they provide owners with a mechanism to describe the actual operating behavior of their bridges. This operating behavior can be used to quantify safety and/or reinforce maintenance decisions such as repair, rehabilitation and replacement alternatives or even a risk-based do nothing alternative.

**Phase I: Intact Element-Level Validation**
In the selected genre of bridges, reinforced concrete slabs and steel girders are the main structural components that comprise in the load-carrying mechanism of the bridge system. As a result, two representative experimental investigations on a corner-supported two-way reinforced concrete slab (20) and a simply-supported steel girder (21, 22) were selected in this study to accomplish the goal of the first phase of the proposed framework. This element-level studies provide a tool for validating the modeling approaches employed through the full nonlinear (material and geometric) response of these components, which provides the confidence for integrating elements into system-level models.

In an experimental study, McNeice (20) performed a scaled test on a 914.4 mm (36 in) square concrete slab having thickness of 44.5 mm (1.75 in). The tested slab was reinforced with an isotropic mesh of steel ($\rho = 0.85\%$) placed in one layer at a depth of 33.3 mm (1.31 in). The slab was supported at four corners and tested under a central load applied at the middle of the span. Due to the symmetry of the structure in geometry and loading conditions, only one quarter of the slab was modeled numerically for the purpose of this study (see Figure 2a). Appropriate boundary conditions were applied at the planes of symmetry as well as a single corner node. In another study by Lagerqvist et al. (21, 22), a series of experimental investigations were conducted on simply-supported high strength steel plate girders under lateral patch loading. Among all tested specimens, a plate girder with geometrical configurations illustrated in Figure 2b, was chosen in this study for numerical model development and validation. Ideal pin-roller boundary conditions were applied at either end of the simulated girder, while the lateral loading was applied linearly across the top flange at the midspan of the girder.

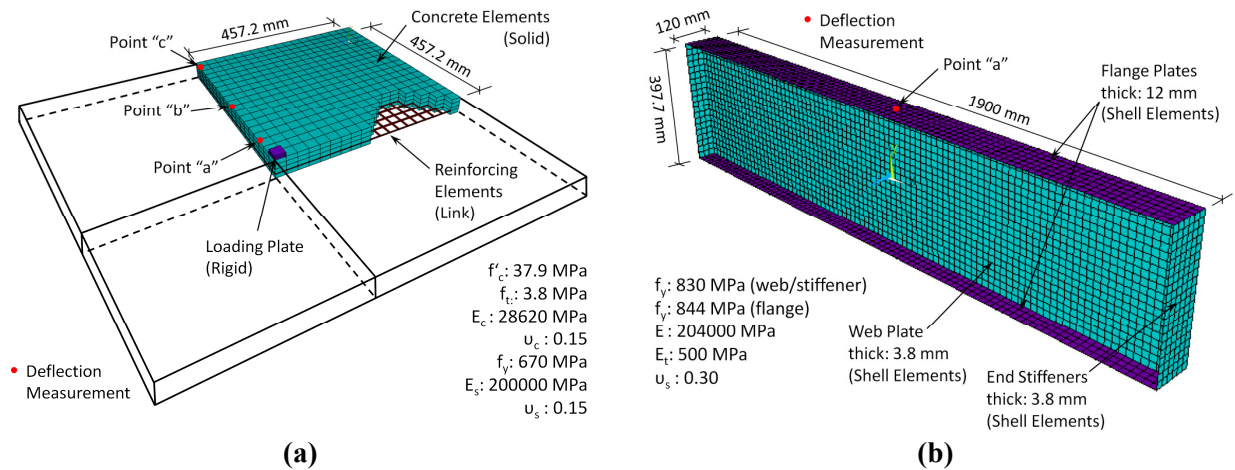

**FIGURE 2 Intact element-level FE model development (a) concrete slab, (b) plate girder.**

Appropriate material constitutive relationships with suitable failure criteria (Von-Mises and William-Warnke (19) for steel and concrete components, respectively) were included in each model to capture the effect of material non-linearities on the element-level behavior. Inelastic stress-strain relationships, cracking and crushing of the concrete, as well as yielding and strain hardening of the steel components are the main sources of material non-linearities incorporated in the models. For the simply-supported steel girder, geometric non-linearity was also included in the analysis to allow for the development of local instabilities, which often dominate the behavior of thin-walled structures (23, 24). Initial imperfections in the format of out-of-plane flatness and twisting were introduced to the web plate and the loaded flange of the simulated plate girder, respectively, to agitate the occurrence of the buckling phenomenon (25).

Displacement-controlled static analysis with the Newton-Raphson non-linear solution algorithm was used for both of the developed models. For the corner-supported concrete slab, load-deflection responses at three predefined nodes (see Figure 2a) were derived from the numerical analysis and compared to the corresponding experimental outcomes. The amount of load applied on the specimen in the experiment was limited to study the serviceability of the slab. However, the numerical analysis was expanded to the post-yield behavior of the steel reinforcement and able to simulate the full non-linear behavior of this element including cracking/crushing in the concrete body and plasticity in steel reinforcement. For the simply-supported steel girder, on the other hand, the experimental study was able to capture the ultimate capacity of the element which was governed by local buckling of the web over most of the member depth (25). The load-deflection behavior at the specific location of the specimen (point "a") was selected for the validation study of the developed numerical model. Figure 3 demonstrates a good correlation between the numerical results and those obtained in the experimental practices.

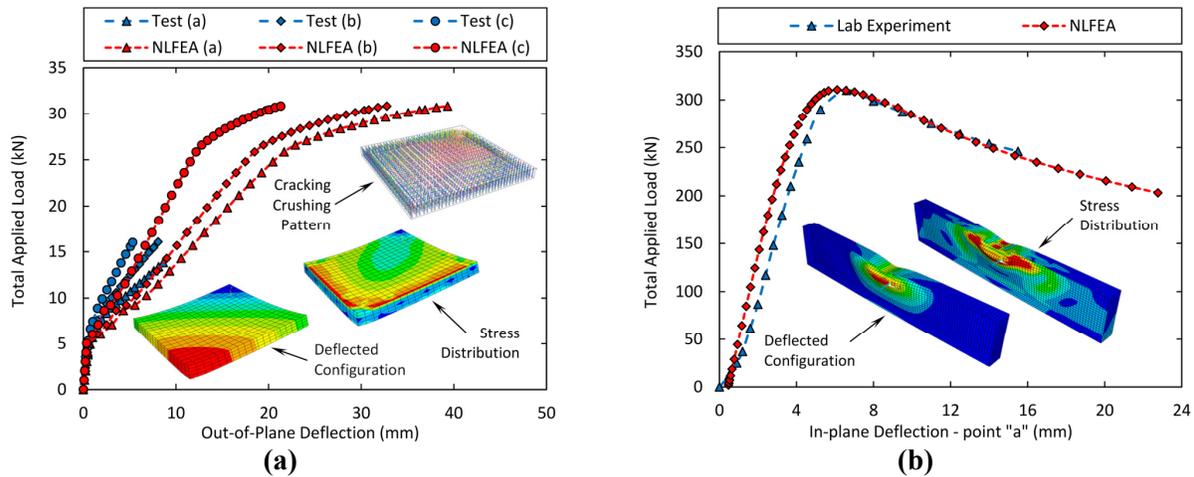
**FIGURE 3 Intact element-level validation (a) concrete slab (b) steel girder.**

**Phase II: Intact System-Level Validation**

A scale laboratory investigation performed at the University of Nebraska (17) was selected in this study with the primary objective of numerically simulating and validating the system-level behavior of composite steel girder bridges. The tested bridge was simply supported with a span length of 21.34 m (70 ft) and cross section width of 7.92 m (26 ft). The superstructure consisted of a 190.5 mm (7.5 in) thick reinforced concrete deck, supported by three steel plate girders spaced at 3.05 m (10 ft), as shown in Figure 4a. An ultimate capacity load test was conducted on this model bridge to evaluate the load carrying capacity of the system. Vertical concentrated loads were applied on the bridge deck to mimic two side-by-side HS-20 trucks. The corresponding loading scenario was longitudinally positioned on the system to produce the maximum positive moment at the mid-span of the bridge.

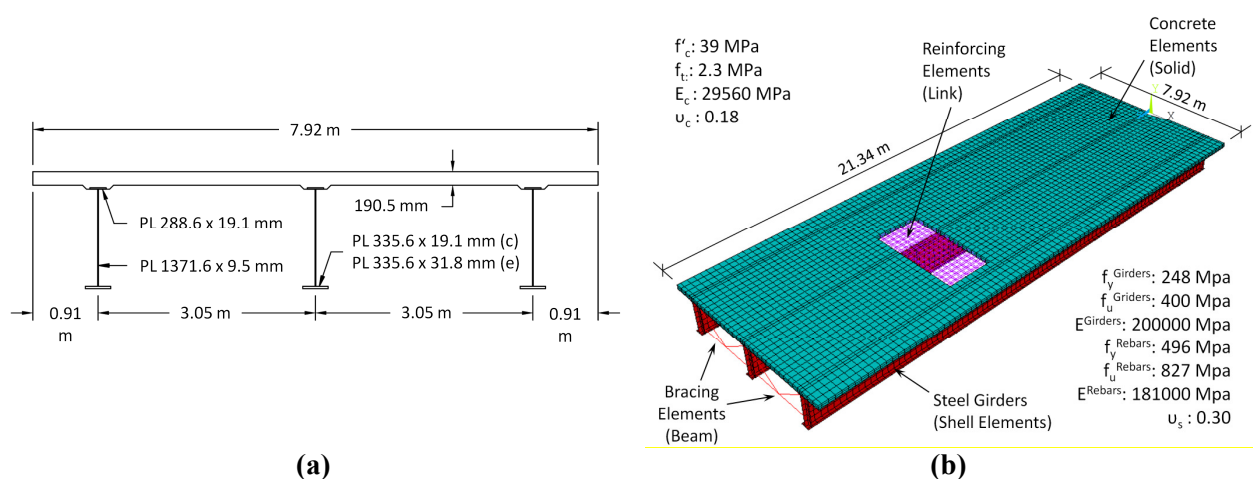
**FIGURE 4 Nebraska laboratory test (a) bridge cross section (b) developed FE model.**

Figure 4b illustrates the FE model for the tested bridge superstructure. All of the structural components including the concrete deck, flexural reinforcement (in two layers), steel girder, and lateral bracing were accurately modeled with the material and section properties given in the test report (17). With no evidence of loss in the composite behavior observed during the test, the girders are assumed to be in full composite action with the concrete deck in the developed model. The model was restrained with hinge and roller supports at the ends and loaded with a series of point loads applied over 500 x 200 mm (20 x 8 in.) patch areas (26). Interior girder deflection at mid-span of the bridge was used to validate the proposed numerical model. As illustrated in Figure 5a, results obtained via non-linear FE analysis correlated well with the experimental outcome. Based on the validated numerical analysis, the general behavior of the investigated bridge was classified into four different stages (see Figure 5b). The bridge system behaves linearly elastic prior to formation of first flexural cracks in the concrete deck (stage A). Beyond first cracking point, it continues to carry more loads until plasticity initiates in steel girders (stage B). With further increment in the external loading, the structural stiffness drops off significantly and plastic hinges are formed in steel girders at the location of maximum moment (stage C). Eventually, the proposed FE model captures the punching shear failure mechanism in the system that terminated the corresponding laboratory experiment (stage D).

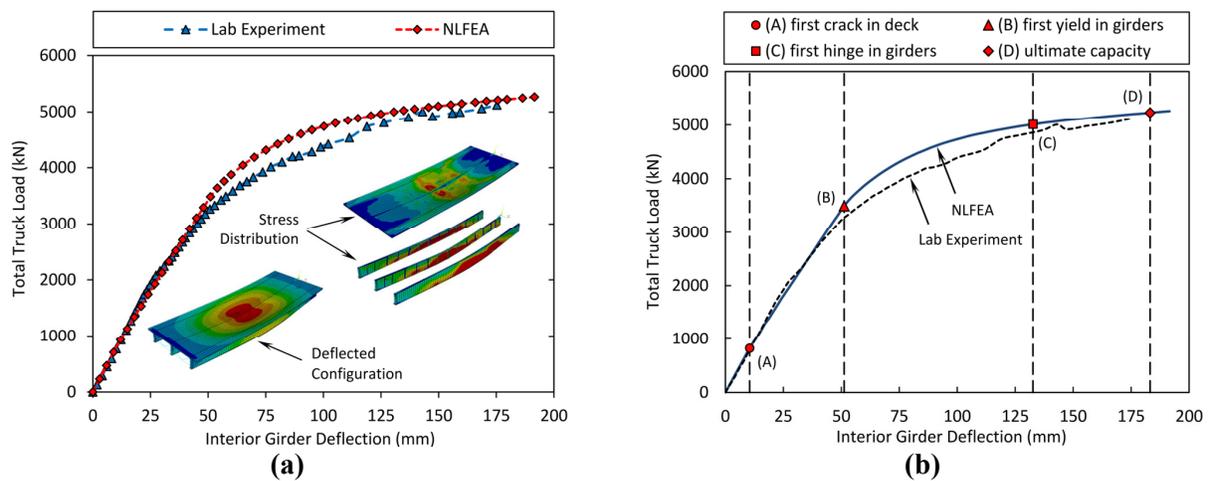

**FIGURE 5 (a) Intact system-level simulation (a) model validation (b) behavioral stages.**

The outcome of this phase was a fundamental understanding of the intricacies necessary for a robust modeling and behavioral characterization of this type of bridge superstructures. The results also provide a rational foundation to characterize the impact of damage and deteriorating conditions on the overall system-level behavior relative to an impact system.

**Phase III: Damaged Element-Level Validation**
With the limited data that exists on the behavior of bridge structures with accumulated damage, the development of modeling strategies for integrating damage at the element-level provides a suitable alternative. Once the element-level damage modeling approach is validated, it can be integrated into the system-level models established in phase II. As a result, the third phase of the proposed framework aims to establish a fundamental understanding to characterize the impact of damage and deteriorating mechanisms on the individual bridge components. In this work, the validation of the damage was limited to girder corrosion, but the process described can be used for other damage mechanisms.

Corrosion in steel girders is one of the common damage scenarios associated with the composite steel girder bridges. This type of damage usually occurs at either end of the girders due to deicing media of salt and water leaking from the deck expansion joints or occasional splashing from passing vehicles. The corresponding deteriorating mechanisms consist of bottom flange thickness reduction, thinning of the bottom portion of the web, and even formation of the irregular-shaped holes in the web just above the bottom flange (27). These deteriorating conditions may significantly reduce the load-carrying capacity of the degraded member and eventually the bridge system, due to localized failure mechanisms such as local buckling, web crippling or crushing of the end stiffeners.

A comprehensive experimental investigation was conducted at Michigan Technological University to characterize the influence of end deterioration on the capacity of steel girders (28). The experimental program was performed on four degraded sub-sections of wide flange beams with variations in the extent of damage. Among these specimens, a W460x158 (W18x106) beam section with the geometrical configurations illustrated in Figure 6a was chosen in this study for model development and validation. This beam has a web damage height of 38.1 mm (1.5 in), flange damage width of 38.1 mm (1.5 in), and a damage depth of 3.2 mm (1/8 in) on each side running all through the length of the specimen. The stiffening plates welded to the top flange and upper portion of the web were meant to be a representative of a deck or other mechanisms which provide rotational restraint and force the buckling to occur at the lower portion of the web where damage exists.

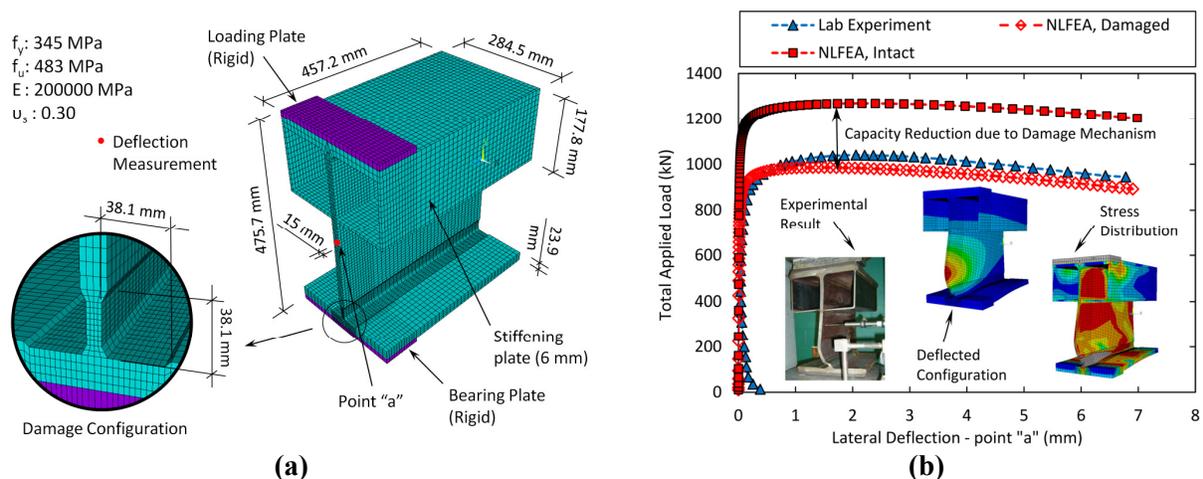

**FIGURE 6 Element-level damage integration (a) model development (b) model validation.**

Two FE models were created to assess the intact and damaged conditions of the selected beam specimen. Both models were fully restrained at the bottom of the bearing plate and uniformly loaded through the top loading plate to mimic the boundary conditions applied in the experiment. Material and geometric non-linearities were both included in the analysis. An initial imperfection compatible with the first buckling mode-shape was also introduced to the models to agitate the buckling phenomenon. Out-of-plane deflection at the mid-height of the section (point "a") was recorded in the test and used in this study to validate the proposed damaged model. As illustrated in Figure 6b, the numerical results of the deteriorated model are well correlated to the experimental outcomes. The modeling approach allowed for simulation of combined global-local buckling failure that was observed in the experiment. Moreover, comparing the numerical results of intact and damaged models demonstrated the magnitude of the influence of the corrosion and section loss damage mechanism on the load-carrying capacity of this element.

**Phase IV: System-Level Damage Integration**
The first three phases aimed to establish a comprehensive foundation on the model development and damage integration strategies, which were essential to the body of the proposed framework.

Accordingly, the main objective of the last phase is to integrate different types of damage and deteriorating mechanisms into the measure of system performance and characterize the impact of damage on the ultimate capacity, redundancy and system ductility. The system-level model with integrated damage would be an extrapolation of the intact model (validated in phase II) coupled with the damage scenario(s) incorporated into the element-level domain (phase III).

To illustrate the impact of damage on the system-level behavior, the bridge model presented in phase II was updated to include the end deterioration mechanism simulated and validated in phase III. As it is depicted in Figure 7a, all three girders were assumed to be deteriorated at both ends with the same level of damage extension. The thicknesses of the web, bottom flange, and lateral stiffeners were uniformly reduced by 40% in the affected areas to mimic a representative corrosion condition that may happen in the actual practice resulting from a deck joint leaking.

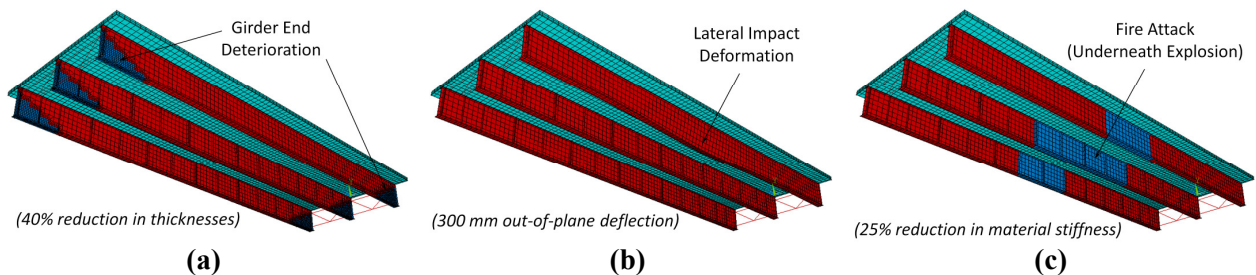

**FIGURE 7 Integrated damage scenarios (a) girder corrosion (b) lateral impact (c) fire.**

In addition to the end girder corrosion, two more damage scenarios were also introduced into the numerical model of the bridge under investigation. The first one was the lateral impact which could be caused by an over height vehicle passing beneath the bridge. This type of collision may result in major out-of-plane deformation in the exterior girder, which was reflected in the model as an initial lateral deflection at the midspan of the bridge (see Figure 7b). The last damage mechanism included in this study was a fire attack that could take place as a result of an explosion of a chemical-carrying truck involved in an accident right underneath the structure. The bridge model was updated to include the effect of fire attack through a 25% reduction in the material stiffness of the girders at the center of the bridge superstructure (see Figure 7c). It should be emphasized that these scenarios are idealized, but provide a foundation for characterizing the system-level influence of likely damage mechanisms.

Under the same loading and boundary conditions applied to the intact system, the updated models were numerically analyzed to evaluate the effect of integrated damage mechanisms on the overall system behavior. Figure 8a illustrates the variations in the system response under the influence of the different damage scenarios. As shown, all the numerical models including the intact system, demonstrated additional reserve capacity over the AASHTO LRFD nominal design capacity (29), which would indicate the high level of inherent redundancy that exists in this type of bridge superstructures. This level of redundancy can be attributed to the complex interaction exists among all structural components of a bridge system which is not considered in the component based analysis approach of current design provisions.

In a given bridge superstructure, whether it is intact or damaged, the measure of system performance can be preliminarily defined based on its ultimate capacity, system redundancy, and system ductility. As it is demonstrated in Figure 8b, the ultimate capacity of a system is the maximum level of load that can be tolerated before the failure mechanism takes place. The difference between the ultimate and design capacity, is characterized herein as system reserve capacity, which would also represent the level of redundancy exists in the structure. Thus, reduction in system reserve capacity due to presence of any damage mechanism can be interpreted as reduction in the system redundancy. Similarly, system ductility can be characterized as the ratio of the maximum measured response of the system at the moment of failure ($\Delta u$) to the system response when the steel girders yield for the very first time ($\Delta y$). However, these definitions are based on the behavioral stages characterized for composite steel girder bridges (Figure 5b) and subjected to change for other types of bridge superstructures.

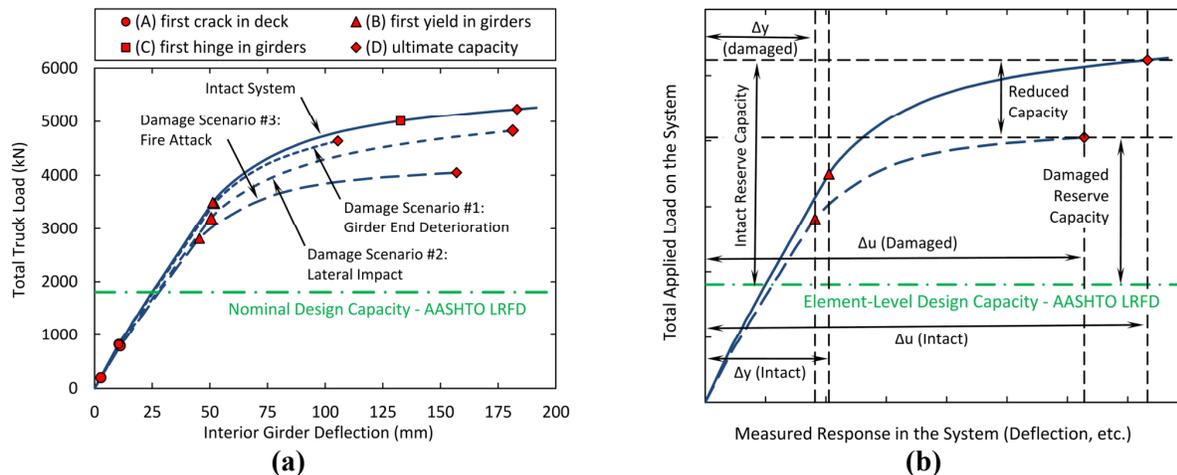

**FIGURE 8 System-level damage integration (a) variations in system behavior (b) measure of system performance.**

Table 1 quantifies the system performance parameters for the analyzed intact and damaged bridge models. The amount of reduction in each parameter of the deteriorated bridges was calculated with respect to the corresponding value of the intact system, (see Figure 8b). As summarized in Table 1, the over-height vehicular impact has the lowest influence on the system redundancy (13.4 %), with negligible effect on the system ductility. The system redundancy was decreased significantly when the bridge superstructure was exposed to the fire (52.9 %), while the end deterioration in the steel girders has the most crucial influence on the system ductility (76.6 %). A reduction in system ductility increases the chance of brittle failure of the structure in the presence of other coupled damage mechanisms or irregular over-sized loading configurations. On the other hand, a reduction in system redundancy would decrease the maximum bearing capacity of the system and eventually diminish the operational safety margin of structure. Depending on the level of reduction, appropriate maintenance decisions shall be made to keep the structure in service. Although the ultimate capacity of the analyzed system decreased in all three cases of incorporated damages, none of the assumed damage scenarios were severe enough to reduce the capacity of system below element-level design limit state.

While these results highlight the robustness in the design of this structure, it should be noted that all calculated parameters are unique to the type, extension and level of the damage mechanisms integrated in the numerical model and would be altered in the presence of a new damage scenario or different damage configuration. In addition, the loading scenarios were not configured to necessarily generate the lower bond on the system response. Possible modifications in truck positions would be shear load configuration in the case of end deterioration and single lane loading over the exterior girder for the impact damage scenario. For a rational comparison of these cases, the behavior of the intact system would also need to be revaluated based on the modified loading configurations.

**Table 1 Measure of system performance.**

| Damage Scenario | Design Capacity (kN) | Ultimate Capacity (kN) | Reserve Capacity (kN) | Reduction in System Redundancy (%) | System Ductility ($\Delta u/\Delta y$) | Reduction in System Ductility (%) |
|---|---|---|---|---|---|---|
| Intact | 1806 | 5235 | 3429 | - | 3.58 | - |
| Girder corrosion | | 4636 | 2829 | 21.2 | 2.02 | 76.7 |
| Lateral impact | | 4829 | 3023 | 13.4 | 3.58 | 0.1 |
| Fire | | 4049 | 2243 | 52.9 | 3.45 | 3.7 |

**FINDINGS AND PATH FORWARD**

The overall objective of this study was to establish a framework to evaluate the in-service condition of bridge superstructures in the presence of common deteriorating mechanisms and provide a measure of system performance by characterizing the impact of damage on the ultimate capacity, redundancy, and system ductility. The damage scenarios included in this study were selected as a representation of deteriorating conditions that may influence the performance and serviceability of highway steel bridges. Thus, all geometrical characteristics of the integrated damage mechanisms including the shape of the affected regions, depth of the damage, and the corresponding extent level were rationally assumed based on literature (10, 27) and integrated into the calibrated bridge model. However, for more comprehensive evaluation of an in-service bridge superstructure, the last phase of the proposed framework would need to be interconnected with a comprehensive non-destructive field inspection for each individual structure selected for evaluation, to accurately model the existing damage condition and its details. The incorporation of condition state data obtained from periodic inspection coupled with the damage-integrated system-level behavior characterization has the potential to provide a real-time estimate of system performance. By updating the developed numerical model of the in-service highway bridges based on biennial inspection data, degradations in the structural performance parameters can be monitored and evaluated over time. Extrapolating the degradation trend through the design life of the structure would help the bridge owners to estimate remaining service life of the bridge system and make appropriate maintenance decisions regarding the long-term preservation strategies.

      This investigation focused exclusively on the composite steel girder bridges and was aimed at representing a conceptual schematic of a computational modeling strategy for describing an in-service baseline performance measure; however, this same methodology could be extended to other types of bridge superstructures. The proposed framework could be beneficial to the preservation community as a mechanism to make decisions based on in-service condition, but also has implications in the design community where a system-level design strategy would have a major impact on design economy as compared to current element-level design strategies. Moreover, results obtained from this investigation highlight the ability of the proposed framework to provide a critical linkage between the design and preservation communities by correlating the element-level behavior to the system-level response under the effect of different damage scenarios. The numerical modeling approach implemented in the proposed framework also has the potential to explore the implication of advances in material, design methodologies and construction practices on the long-term performance of bridge superstructures. All things considered, the future path of this research study, can be gravitated

towards implementation of the proposed framework on a series of representative existing bridges with a comprehensive inspection reports to evaluate their in-service performance and help owners to make appropriate maintenance decisions.